\documentclass[11pt]{article}
\usepackage{amsmath,amsfonts,amssymb,graphicx,indentfirst}
\usepackage{natbib}
\bibpunct{(}{)}{;}{a}{,}{,}
\newcommand{\eps}{\epsilon}
\newcommand{\ups}{\upsilon}

\newcommand{\Exp}[2][]{\mathbf{E}_{#1}\left(#2\right)}

\newcommand{\definition}{\overset{\text{def}}{\Leftrightarrow}}

\newcommand{\tendsto}{\rightarrow}
\newcommand{\me}{\mathrm{e}}

\newcommand{\sN}{\sigma_N}
\newcommand{\mN}{\mu_N}
\newcommand{\kN}{\kappa_N}

\newcommand{\Ai}{\mathrm{Ai}}
\newcommand{\Bi}{\mathrm{Bi}}
\newcommand{\AB}{\mathbf{AB}}
\newcommand{\id}{\mathrm{Id}}
\newcommand{\wkl}{W_{\kappa,\lambda}}

\newcommand{\mE}{\mathbf{E}}
\newcommand{\mM}{\mathbf{M}}
\newcommand{\cE}{{\cal E}}
\newcommand{\cM}{{\cal M}}

\newcommand{\mEi}{\mathbf{E}^{-1}}
\newcommand{\cEi}{{\cal E}^{-1}}

\newcommand{\gO}{\mathrm{O}}
\newcommand{\lo}{\mathrm{o}}

\newcommand{\lb}{\lambda\beta^2}
\newcommand{\pcf}{parabolic cylinder functions}
\newtheorem{theorem}{Theorem}
\newtheorem{fact}{Fact}
\title{On the largest eigenvalue of Wishart matrices with identity covariance when $n$, $p$ and $p/n \tendsto \infty$}
\author{Noureddine El Karoui\thanks{\textbf{Acknowledgements:} The author is grateful to Pr. Iain Johnstone for
 many discussions, key references, and guidance and to Pr. David Donoho for his helpful comments and support.
Supported in part by NSF \textsf{DMS}-0140698 and \textsf{ANI}-008584 (ITR).  \textbf{AMS 2000 SC: } Primary
62E20, Secondary 62H25. \textbf{Key words and Phrases : } Principal Component Analysis, largest singular value,
Tracy-Widom distribution, Fredholm determinant, Random Matrix Theory, Wishart Matrices. \textbf{Contact :} \texttt{nkaroui@stanford.edu}} \\
\textit{Department of Statistics,}\\ \textit{ Stanford University} }

\begin{document}
\maketitle
\begin{abstract}
Let $X$ be a $n\times p$ matrix and
$l_1$ the largest eigenvalue of the covariance matrix
$X^{*}X$. The ``null case" where $X_{i,j}\sim {\cal N}(0,1)$ is of particular interest for principal component analysis.

For this model, when $n, p\tendsto \infty$ and $n/p \tendsto \gamma \in \mathbb{R}_+^*$, it was shown in
\citet{imj} that $l_1$, properly centered and scaled, converges to the Tracy-Widom law.

We show that with the same centering and scaling, the result is true even when $p/n$ or $n/p\tendsto \infty$,
therefore extending the previous result to $\gamma \in \overline{\mathbb{R}}_+$. The derivation uses ideas and
techniques quite similar to the ones presented in \citet{imj}. Following \citet{sosh}, we also show that the same
is true for the joint distribution of the $k$ largest eigenvalues, where $k$ is a fixed integer.

Numerical experiments illustrate the fact that the Tracy-Widom approximation is reasonable even when one of the
dimension is small.
\end{abstract}
\section{Introduction}
Large scale principal component analysis (PCA) - concerning an $n\times p$ matrix $X$ where $n$ and $p$ are both
large - is nowadays a widely used tools in many fields, such as image analysis, signal processing, functional data
analysis and quantitative finance. Several examples come to mind, including Eigenfaces, subspace filtering, or
\citet{lalouxetal} where PCA (as well as some random matrix theory) is used to try to improve on the naive
solution to Markovitz's portfolio optimization problem.

Important progress has been made recently in our understanding of the statistical properties of PCA in such
settings. Emblematic of this is work of \citet{imj}, which explains the properties of the square of the largest
singular value of a random matrix $X$ under the ``null model" where its entries are iid
${\cal N}(0,1)$. Specifically, if we denote the sample eigenvalues of $X'X$ by $l_1\geq \ldots \geq l_p$,   call
\begin{align*}
n_1&=\max{(n,p)}-1\;, \; \; \; \; p_1=\min{(n,p)} \;,\\
\mu_{np}&=(\sqrt{n_1}+\sqrt{p_1})^2 \; ,\\
\sigma_{np}&=(\sqrt{n_1}+\sqrt{p_1})\left(\frac{1}{\sqrt{n_1}}+\frac{1}{\sqrt{p_1}}\right)^{1/3} \; ,
\end{align*}
and $W_1$ the Tracy-Widom distribution (see \textbf{A0}), it was shown in \citet{imj} that
\begin{theorem}[Johnstone]\label{ThJ}
If $n,p \tendsto \infty$ and $n/p\tendsto \gamma \in (0,\infty)$,
$$
\frac{l_1-\mu_{np}}{\sigma_{np}} \overset{\cal L}\tendsto W_1 \;.
$$
\end{theorem}

Building on \citet{imj} and using properties of determinantal point processes, \citet{sosh} showed that the same result
holds for the $k$ largest eigenvalues, where $k$ is a fixed integer: their joint distribution converges to their
Tracy-Widom counterpart.

This is a very interesting development because the classical theory (e.g \citet{anderson}) was developed under the
assumption that
$p$ was fixed and $n$ grew to $\infty$, whereas in modern day applications both $p$ and $n$ are large. However,
Johnstone's assumption $n/p\tendsto \gamma$ imposes a limit on the validity of his result which one would like to
remove. In an actual data analysis, with given $p$ and $n$, $n=\lo(p)$ and $n\asymp p$ could be equally plausible.
Furthermore, a specific $X$ of size $n\times p$ could arise in many triangular arrays settings, where we have
$X_j$ of size $n_j\times p_j$, and the limitation
$n_j/p_j\tendsto \gamma$ finite might only hold in some triangular situations and not in others.

Accordingly in this paper we weaken the assumption that $n/p \tendsto \gamma$ finite and show that
\begin{theorem}\label{ThJextended}
If $n,p \tendsto \infty$ and $n/p\tendsto \infty$,
$$
\frac{l_1-\mu_{np}}{\sigma_{np}} \overset{\cal L}\tendsto W_1 \;.
$$

Moreover, with the same centering and scaling, the joint distribution of the $k$ largest eigenvalues converges in
law to its Tracy-Widom counterpart.

 Dually, the same result holds if $n/p\tendsto 0$.
\end{theorem}
Let us note that the remark we made about centering and scaling sequences after Theorem~\ref{ThJ} is still valid in
this context.

There is clearly a mathematical motivation for dealing with this problem: the result completes the picture about the
properties of $l_1$ with large $p$ and $n$ and, in a sense, closes Theorem \ref{ThJ}. But is it interesting from a
statistical standpoint?

The situation $p\gg n$ is indeed a fairly common one in modern statistics. Microarray data are a prototypical example:
currently they usually have $p$ of the order of a few thousands and $n$ of the order of a few tens. One encounters
$p\gg n$ or $n\gg p$ in many other instances: data collection mechanisms are now effective enough so as to, for
example, collect and retain thousands of piece of information for millions of customers (transactional data), or
millions of pieces of information for thousands of stocks (tick-by-tick data in Finance). Analyzing these very high
dimensional datasets raises new challenges and is at the center of recent statistical work, both applied and
theoretical.

Microarray analysis in particular is a very active field, and has contributed a flurry of activity in non
classical situations (very high dimensional data), raising theoretical questions and sometimes revisiting
classical techniques or results. As illustrated for instance in \citet{wrr}, PCA or PCA-related methods are used
for various tasks in the microarray context, from traditional dimensionality reduction procedures to gene
grouping. Having a good understanding of the behavior of the singular values  of gaussian ``white noise" matrices
could provide valuable insights for these applications. Recent work of \citet{bl03} about the properties of naive
Bayes and Fisher's linear discriminant function when $p\gg n$ illustrates the impetus these dimensionality
assumptions are also gaining in theoretical studies. Our work is part of the larger effort to investigate the
properties of high dimensional data structures. Here it is done in a simple, ``null" situation.

We now present a few numerical experiments we realized to assess how big (or small) $n$ or $p$ should be for
Theorems \ref{ThJ} and \ref{ThJextended} to be practically useful.

\subsection{Numerical experiments}
\citet{imj} showed empirically that in that situation the Tracy-Widom approximation was reasonably satisfying, even for
small matrices. Similarly, to try to assess its accuracy in our setup, we ran the following experiments in
\texttt{Matlab}: we picked $n$ and $p$ and generated $10,000$ $n\times p$ matrices $X$ with entries iid ${\cal
N}(0,1)$. Then we used standard routines (\texttt{normest} in \texttt{Matlab}) to compute their spectral norms and
squared them to obtain a dataset of $l_1$-s.

Following \cite{imj2}, we adjust centering and scaling to
\begin{align*}
\tilde{\mu}_{np}&=\sqrt{n-1/2}+\sqrt{p-1/2}\; ,\\
\tilde{\sigma}_{np}&=(\sqrt{n-1/2}+\sqrt{p-1/2})\left(\frac{1}{\sqrt{n-1/2}}+\frac{1}{\sqrt{p-1/2}}\right)^{1/3}
\;.
\end{align*}
This leads to a very significant improvement in the quality of the Tracy-Widom approximation for our simulations.
Simple manipulations (explained in section 2.2) show that we have some freedom in choosing the centering and
scaling: if we replace $n$ by
$n+a$ and $p$ by $p+b$ (where $a$ and $b$ are fixed real numbers) in the definitions of $\mu_{np}$ and
$\sigma_{np}$, Theorem~\ref{ThJ} and Theorem~\ref{ThJextended} still hold.
The particular choice used here is motivated by a careful theoretical analysis of the entries of $K_N$ mentioned
in section 2.2.

%Remark that $n$ and $p$ now play completely symmetric roles, which was not the case with $\mu_{np}$ and
%$\sigma_{np}$: there, the larger dimension was not treated as the .

Table 1 summarizes the ``quantile" properties of the empirical distributions we obtained and compare them to the
Tracy-Widom reference. We used the same reference points as \citet{imj}.

\begin{table}
\begin{center}
\begin{tabular}{|c|c|ccc|cc|}\hline
TW Quantiles & TW & 10$\times$1000 & 10$\times$ 4000 &10$\times$ 10000&100$\times$4000 &30$\times$5000
\\ \hline
-3.9 &.01 & 0.009  &  0.010  &  0.015  &  0.012  &  0.013  \\
 -3.18 &.05 &0.047  &  0.050  &  0.060  &  0.053  &  0.055  \\
-2.78 &.10 &    0.102  &  0.107  &  0.112  &  0.103  &  0.105  \\
-1.91 & .30 &  0.303  &  0.308  &  0.316  &  0.304  &  0.303  \\
-1.27 & .50 & 0.506  & 0.506   &  0.522  &  0.508  &  0.503  \\
-0.59 & .70 & 0.705  & 0.704 &  0.723  & 0.706  & 0.702  \\
0.45 & 0.9 & 0.904  & 0.904 &  0.913  & 0.901  & 0.904  \\
0.98 & .95 & 0.953  & 0.951 &  0.958  & 0.951  & 0.953  \\
2.02 & .99 & 0.992  & 0.990 &  0.992  & 0.991  & 0.991  \\\hline

\end{tabular}
\\[.5cm]
\begin{tabular}{|c|c|ccc|ccc|}\hline
TW Quantiles & TW & 50$\times$5000 & 50$\times$20000 & 50$\times$50000 & 5$\times$200 &5$\times$2000& 5$\times$20000\\
\hline

-3.9 &.01 & 0.010 &   0.017 &   0.021 &   0.008 &   0.014 &   0.018 \\
-3.18 &.05 & 0.053 &   0.067 &   0.079 &   0.047 &   0.057 &   0.069 \\
-2.78 &.10 &  0.104 &   0.125 &   0.139 &   0.094 &   0.110 &   0.120\\
-1.91 & .30 &  0.309 &   0.331 &   0.345 &   0.293 &   0.314 &   0.320\\
-1.27& .50 &   0.502 &   0.522 &   0.538 &   0.500 &   0.506 &   0.519\\
-0.59 & .70 &  0.705 &   0.718 &   0.727 &   0.714 &   0.712 &   0.710\\
0.45 & .90 &   0.899 &   0.905 &   0.911 &   0.911 &   0.906 &   0.907\\
0.98 & .95 &   0.949 &   0.955 &   0.957 &   0.959 &   0.951 &   0.954 \\
2.02 & .99 &   0.991 &   0.992 &   0.992 &   0.994 &   0.992 &  0.992\\ \hline
\end{tabular}
\end{center} \caption{\textbf{Quality of the Tracy-Widom Approximation for some large matrices:} the leftmost columns
displays  certain quantiles of the Tracy-Widom distribution. The second column gives the corresponding value of
its cdf. Other columns give the value of the empirical distribution functions obtained from simulations at these
quantiles. $\tilde{\mu}_{np}$ and $\tilde{\sigma}_{np}$ are the centering and scaling sequences.}
\end{table}
We picked the dimensions according to two criteria: $100\times4000$, $30\times5000$, and $50\times5000$ were
chosen to investigate ``representative" microarray situations. We chose the other to have a range of ratios and
estimate how valuable the Tracy-Widom approximation would be in situations that could be considered classical, i.e
one small dimension (less than 10) and one large (several hundreds to several thousands). For the sake of
completeness, we redid the simulations presented in \citet{imj} and present in Table 2 the results obtained with
$\tilde{\mu}_{np}$ and $\tilde{\sigma}_{np}$ as centering and scaling.

\begin{table}
\begin{center}
\begin{tabular}{|c|c|ccc|ccc|}\hline
TW Quantiles & TW & 5$\times$5 & 10$\times$ 10 &100$\times$ 100&5$\times$20 &10$\times$ 40 & 100$\times$400
\\ \hline
-3.9 &.01 &0      & 0.002   &  0.008  &  0.001  &  0.004 &  0.008 \\
-3.18 &.05 & 0.003  & 0.018   &  0.043  &  0.019  &  0.032 & 0.044 \\
-2.78 &0.10 &0.022  &  0.054  &  0.090  &  0.056  &  0.077 &  0.095 \\
-1.91& .30 & 0.217  &  0.257  &  0.295  &  0.262  &  0.279 &   0.294 \\
-1.27 & .50 & 0.464  &  0.486  &  0.497  &  0.490  &  0.494 &  0.489 \\
-.59 & .70 & 0.702  &  0.703  &  0.700  &  0.702  &  0.707 &   0.702 \\
0.45 & .90 & 0.903  &  0.903  &  0.901  &  0.905  &  0.906 &   0.899 \\
0.98 & .95 & 0.949  &  0.950  &  0.950 &  0.952  &  0.953 &  0.949 \\
2.02 & .99 & 0.988  &  0.990  &  0.991  &  0.989  &  0.990 &   0.990 \\ \hline
\end{tabular}
\end{center}
\caption{\textbf{Quality of the Tracy-Widom Approximation (Continued):} the columns have the same meaning as in
Table 1. The ratio  $p/n$ is smaller than in Table 1 and the matrices are not as big, but the Tracy-Widom
approximation is already acceptable for the upper quantiles. $\tilde{\mu}_{np}$ and
$\tilde{\sigma}_{np}$ are the centering and scaling sequences. }
\end{table}

We see that the fit is good to very good for the upper quantiles ($.9$ and beyond) across the range of dimensions
we investigated. The practical interest of this remark is clear: these are the quantiles one would naturally use
in a testing problem. We note that it appears empirically that the problem gets harder when the ratio
$r$ of the larger dimension to the smaller one ($p_1$ in our notation) gets bigger: the larger $r$, the larger $p_1$
should be for the approximation to be acceptable.

\subsection{Conclusions and Organization}
From a technical standpoint, the method developed in \citet{imj} proves to be versatile, and, at least
conceptually, relatively easy to adapt to the case where $n/p\tendsto \infty$. Nevertheless, substantial technical
work is needed to obtain Theorem~\ref{ThJextended}. Using the elementary fact (see e.g theorem 7.3.7 in
\citet{hj}) that the largest eigenvalue of
$X^{*}X$ is the same as the largest eigenvalue of
$XX^{*}$, it will be sufficient to give the proof in the case $n/p\tendsto \infty$.

 From a practical point of view, we show that the Tracy-Widom limit law does not
depend of how the sequence
$(n,p)$ is embedded. As long as both dimensions go to infinity, the properly re-centered and re-scaled largest
eigenvalue converges weakly to this law. \\
We can compare this with the ``classical" situation where
$p$ is held fixed, in which case the limiting joint distribution is known, too (see e.g \citet{anderson}, corollary
13.3.2). In this case, the centering is done around $n$ and the scaling is $\sqrt{n}$; elementary computations
show that $(l_1-\mu_{np})/\sigma_{np}$ also has a non-degenerate limiting distribution (possibly changing with
each $p$). Nevertheless, even with the classical centering, it is hard to evaluate the marginals in this context
and the results are therefore difficult to use in practice.

Our simulation results show that the Tracy-Widom approximation is reasonably good (for the upper quantiles) even
when $p$ or $n$ are small. As remarked by \citet{imj}, Proposition 1.2, this implies that when doing PCA, one
could develop (conservative)
tests based on the Tracy-Widom distribution that could serve as alternatives to the scree plot or the Wachter plot.\\

The paper is organized as follows: after presenting (Section 2) the main elements of the proof of
Theorem~\ref{ThJ}, we describe (Section 3) the strategy that will lead to the proof of Theorem~\ref{ThJextended}.
We prove the two crucial points needed in Section 4. To make the paper self-contained, we give some background
information about different aspects of the problem in the appendices. Several technical issues are also treated
there in order to avoid obscuring the proof of the main result.

\section{Outline of Johnstone's proof}
Before describing the backbone of the proof presented in \citet{imj}, we need to introduce a few notational
conventions. In what follows, we will use $N$ instead of $p$ to be consistent with the literature. We also denote
by $\AB$ (for ``asymptotic behavior") the situation where $n,N, \text{ and } n/N \tendsto \infty$. We will
frequently index functions that depend on both $N$ and $n$ with only $N$. The reason for this is that it will
allow us to keep the notations relatively light, and that we think of $n$ as being a function of $N$. Notations
like
$\textbf{E}_N$ and $\textbf{P}_N$ will denote expectation and probability under the measure induced by the
matrices
(of size $n(N)\times N$) we are working with.\\
Finally, it is technically simpler to work with a matrix $X$ whose entries are standard complex Gaussians (i.e the
real and imaginary parts are independent, and they are both ${\cal N}(0,1/2)$), rather than with entries that are
${\cal N}(0,1)$. When we mention the complex case, we refer to this situation.

We now give a quick overview of the important points around which the proof of Theorem \ref{ThJ} was articulated.

At the core of several random matrix theory results lie the fact that the joint distribution of the eigenvalues of
the random matrices of interest is known and can be represented as the Fredholm determinant of a certain operator
(or a totally explicit function of it).

Building on this, if we introduce a number $b$ that is $1$ in the real case and $2$ in the complex one, it turns
out that one has the representation formula
\begin{equation}\label{freddet}
\Exp[N]{\prod_{i=1}^N (1+f(l_i))}=\left[\det(\id+S_Nf)\right]^{b/2} \; ,
\end{equation}
where $S_N$ is an explicit kernel, depending of course upon the kind of matrices in which one is interested. Here,
$f$ treated as an operator means multiplication by this function. It is clear that if $\chi_t=-\mathbf{1}\{x:x\geq t\}$, we
have
$$
\mathbf{P}_N(l_1\leq t)=\left[\det(\id+S_N\chi_t)\right]^{b/2}\;.
$$

The interested reader can find background information on this in \citet{mehta}, chapters 5 and 6, \citet{tw98} or
\citet{deift}, chapter 5, which in turn (p.109) points to \citet{simonreed}, section 17, vol 4, for background on
operator determinants. We stress the fact that all these formulas are finite dimensional.

From the last display, the strategy to show convergence in law in either Theorem \ref{ThJ} or \ref{ThJextended} is
clear: fix $s_0$, show that under the relevant assumptions, $\mathbf{P}(l_{1,N}\leq s_0)\tendsto W_1(s_0)$, and
use the fact that
$W_1$ is continuous to conclude.
\subsection{Complex case}
%In the case we are interested in, the representation (\ref{freddet}) will be exhibited in detail in the complex
%case (for the sake of technical convenience).
We just saw that to find the asymptotic behavior of $l_1$ is equivalent to showing the convergence of the
determinant of a certain operator. This task can be reduced to showing convergence in trace class norm of this
operator (see \citet{simonreed} for background on this, e.g, Lemma XIII.17.4 (p.323)). Through work from
\citet{widom99}, \citet{imj} exhibits an integral representation formula for his operator, and the original
problem is essentially transformed into showing that certain integrals have a predetermined limit.
%The approach
%for this last step is standard: show pointwise convergence and dominate by a function that converges quickly to
%zero.

In somewhat more detail, if we call $\alpha=n-N$, and $L_k^{\alpha}$ the $k$-th Laguerre polynomial associated
with
$\alpha$ (as in \citet{szego}, p.100), let
$$
\phi_k(x)=\sqrt{\frac{k!}{(k+\alpha)!}}x^{\alpha/2}\me^{-x/2}L_k^{\alpha}(x) \; ,
$$
$\xi_k(x)=\phi_k(x)/x$, $a_N=\sqrt{Nn}$, and finally
$$
\left\{
\begin{array}{cll}
\phi(x)&=(-1)^N\sqrt{\frac{a_N}{2}} (\sqrt{n}\xi_N(x)-\sqrt{N}\xi_{N-1}(x))  \; ,\\
\psi(x)&=(-1)^N\sqrt{\frac{a_N}{2}}(\sqrt{N}\xi_N(x)-\sqrt{n}\xi_{N-1}(x)) \; .
\end{array}
\right.
$$
We note two things: first, there is a slight abuse of notation since $\phi$ and $\psi$ obviously depend on $n$ and
$N$, but as in \citet{imj}, we choose to not carry these indices in the interest of readability. Also, $\phi$ and
$\psi$ admit more ``compact" representations, in terms of a single Laguerre polynomial, with a modified $\alpha$, or
another degree. These are easy to derive using \citet{szego}, p.102, for instance. Nevertheless we choose to work
(except in \textbf{A7}) with the previous representations because of the symmetries they present.

The kernel $S_N$ mentioned in (\ref{freddet}) has the representation (\citet{imj}, equation (3.6))
$$
S_N(x,y)=\int_{0}^{\infty} \phi(x+z)\psi(y+z)+ \psi(x+z)\phi(y+z) dz \; .
$$

Now let $\bar{S}$ be the Airy operator. Its kernel is
$$
\bar{S}(x,y)=\frac{\Ai(x)\Ai'(y)-\Ai(y)\Ai'(x)}{x-y}=\int_0^{\infty}\Ai(x+u)\Ai(y+u) du \; ,
$$
where Ai denotes the Airy function. It was shown in \citet{tw94} that, viewing
$\bar{S}$ as an operator on $L^{2}[s,\infty)$, one had
$$
\det(\id-\bar{S})=W_2(s) \; ,
$$
where $W_2$ is the Tracy-Widom law ``emerging" in the complex case (see \textbf{A0}). So the complex analog of
theorem \ref{ThJ} follows from the fact that, after defining $S_{\tau}(x,y)=\sN S_N(\mu_N+\sN x,\mu_N+\sN y)$,
Johnstone managed to show, for all
$s$, that
$$
\det(\id-S_{\tau}) \tendsto \det(\id-\bar{S}) \; .
$$

 To do this, he introduced
 $\phi_{\tau}(s)=\sigma_N \phi(\mu_N+s\sigma_N)$, and
similarly $\psi_{\tau}$. Note that we have
$$
S_{\tau}(x,y)=\int_{0}^{\infty} \phi_{\tau}(x+z)\psi_{\tau}(y+z)+ \psi_{\tau}(x+z)\phi_{\tau}(y+z) dz \; .
$$
Since what we are interested in is really $S_{\tau} \chi_s$, for some fixed $s$, we will view $S_{\tau}$ as an
operator
acting on $L^2[s,\infty)$ in what follows. \\
So the problem becomes to show that, as $n,N \tendsto \infty$
\begin{equation} \label{whattoshow1}
\phi_{\tau}(s),\psi_{\tau}(s)\tendsto \frac{1}{\sqrt{2}}\mathrm{Ai}(s) \; ,
\end{equation}
 and that $\forall s_0 \in \mathbb{R}$,  there exists $N_0(s_0)$ such that if $N>N_0$, we have  on $[s_0,\infty)$,
\begin{equation} \label{whattoshow2}
\phi_{\tau}(s),\psi_{\tau}(s)=O(\me^{-s/2}) \; .
\end{equation}
Once this is shown (we give more details on this later), we can show that $S_{\tau} \tendsto \bar{S}$ in the trace
class norm of operators on $L^2[s,\infty)$. A classical way to do it is described in the remark at the end of
section 3 of \citet{imj}, which bounds the trace class norm of the difference of
$S_{\tau}-\bar{S}$ in terms of the Hilbert-Schmidt norm of operators whose kernels are related to $\phi_{\tau},
\psi_{\tau}$ and
$\Ai$.  This leads to the conclusion that
$$
\det(\id-S_{\tau}) \tendsto \det(\id-\bar{S}) \;,
$$
since $\det$ is continuous with respect to trace class norm. Therefore, the largest eigenvalue of $X^*X$ has the
behavior it was claimed it has.
\subsection{Real Case}
In the real case, using arguments from \citet{tw96} and \citet{widom99}, \citet{imj} gets a representation similar
to (\ref{freddet}), this time involving an operator with kernel a $2\times 2$ matrix (instead of scalar in the
complex case). He is then able to relate it to the complex case problem  - the matrix operator determinant can be
computed as the product of two scalar operator determinants - and shows that the ``reduced" variable he works with
ought to have the same limit as it had in the Gaussian Orthogonal Ensemble case, which was studied in depth by
Tracy and Widom.

For the sake of completeness, we recall that in this situation $\alpha=n-1-N$ and
$$
\mathbf{P}_N(l_1\leq t)=\sqrt{\det(\id+K_N\chi_t)} \;.
$$
$K_N$ has the representation (in the $N$ even
case)
$$
K_N=
\begin{pmatrix}
S_N+ \psi\otimes \epsilon \phi & S_N D-\psi\otimes\phi \\
\eps S_N -\eps +\eps \psi \otimes \eps \phi & S_N + \eps \phi \otimes \psi
\end{pmatrix} \; ,
$$
where $D$ is the differential operator, $\eps$ is convolution with the kernel $\eps(x-y)$, and
$\eps(x)=\text{sgn}(x)/2$. We note the slight change in $\alpha$ and replace $n$ by
$n-1$ when we need to use the results or formulas derived in the complex case
(for instance, the $S_N$ we just mentioned is $S_{n-1,N}$, and not $S_{n,N}$). We refer the reader to \citet{ggk}
for a complement of information on operator determinants and to the end of section VIII in \citet{tw96} for
details on the technical problems that $K_N$ poses.

From a purely technical standpoint, one critical issue is to evaluate the large
$n,N$ limit of
$c_{\phi}=\int_{0}^{\infty}\phi(x)dx/2$. If one can show that it is $1/\sqrt{2}$ when $N\tendsto \infty$ through even
values, then Johnstone's considerations hold true all the way and we have the
same conclusion as in Theorem \ref{ThJ}.\\
We note that using the interlacing properties of the singular values (as mentioned for instance in \citet{sosh},
Remark 5; see also \citet{hj}, theorem 7.3.9), as well as the estimates of the difference (resp. ratio) between
two consecutive terms of the centering (resp. scaling) sequence, the $N$ odd case follows immediately from the $N$
even case. To be more precise, we use the fact that
$$
\frac{\mu_{n,N}-\mu_{n,N-1}}{\sN}=\gO(N^{-1/3}) \; \; \tendsto 0 \text{ as } N \tendsto \infty
$$
to check that the $N$ even terms lower and upper bounding the $N$ odd probability have the same limit. Note that
the same relationship holds for $\mu_{n+a,N+b}$ and $\mu_{n,N}$, if $a$ and $b$ are fixed real numbers. Therefore,
after doing the proof with centering sequence $\mu_{n+3/2,N+1/2}$ (which is technically simpler), we will be able
to conclude that the theorem holds true for
$\mu_{n,N}$.

Last, to be able to use \citet{sosh}, Lemma 2, which gives the result we wish for the joint distribution of the
$k$-largest eigenvalues, we will need to verify that the entries of the $2\times 2$ operator converge pointwise, and
are bounded above in an exponential way. This is what is done in the proof of Lemma 1 of \citet{sosh}, and we will
show in \textbf{A8} that the arguments given there can be extended to handle our situation.

\section{Further Remarks and Agenda}
Most of the work in \citet{imj} is done in closed form, and in the finite dimensional case. That has two
advantages from our standpoint: as the limiting behavior is only investigated in the last ``step", most of the
arguments given there carry through for our problem, and the method certainly does.

 Therefore, our contribution is  mostly technical;
it follows very closely the ideas of \citet{imj}, providing solutions to technical problems appearing in the case
we consider. Only at a few points could we not use the approach developed in \citet{imj}. This led us to an
analysis of the complex case that is slightly different from the original one, but the core reasons for which the
result holds are the same.

In what follows, we first focus on showing that (\ref{whattoshow1}) and (\ref{whattoshow2}) hold true when $n,N$
and their ratio tend to infinity. This takes care of the complex case. We then turn to the problem of the
asymptotic behavior of
$c_{\phi}$, and the technical points we have to verify for \citet{sosh} results to hold.

The following remarks outline the differences between the analysis we present here and the one done in
\citet{imj}.
\subsection{Remarks on adaptation of the original proof}
\subsubsection{Complex case}\label{errorcontrolproblem}
To show that (\ref{whattoshow1}) and (\ref{whattoshow2}) held true, \citet{imj} essentially reduced his problem to
studying the solution of a ``perturbed" Airy equation and used tools from \citet{olver} to carefully study it. One
point that was used repeatedly was that the turning points of the equation were bounded away from one another when
$n,N$ were large. This is not true anymore in the case we consider, and we show how to get around
this difficulty. So we do not work with a perturbed Airy equation anymore, but rather with Whittaker functions,
which have a close relationship to Laguerre polynomials, and their expansion in terms of parabolic cylinder
functions (see \textbf{A9} for some background information on special functions). In \citet{olver80}, the case we
are interested in was studied in detail, giving us most of the tools we need to show (\ref{whattoshow2}). Using
\citet{olver75}, we reinterpret the \pcf\ results in terms of Airy functions and derive the elements we need to
complete the proof of (\ref{whattoshow1}) and (\ref{whattoshow2}).

The reason for which we could not exactly follow the ``original" method is related to the error control function
called
${\cal V}(\zeta)$ in \citet{imj}. This function depends upon the parameter $\omega=2\lambda/\kappa$, which in the case
$n/N\tendsto \gamma \in \mathbb{R}$ is bounded away from 2. This essentially allows a uniform control over ${\cal V}$,
and it is possible to show that this error control function is bounded as a function of $N$. Since the control is
actually something like $\exp(\lambda_0 {\cal V}/\kappa)-1$, it tends to
zero as $N\tendsto \infty$. This gave \citet{imj} a way to get part of (\ref{whattoshow2}). \\
In our case, it seems that ${\cal V}$ would tend to $\infty$, at a rate that is nevertheless $\lo(\kappa)$. As it
seems easier and more promising to use \citet{olver80} than to derive the growth of ${\cal V}$, we choose this
approach. Nevertheless, this is the only (but crucial) technicality (in the complex case) that did not carry
through by the method described in \citet{imj} under $\AB$.
\subsubsection{Real Case}
For the $c_{\phi}$ problem, we provide a closed form expression at given $n,N$ and show that in the limit is the
``right" one as long as $n$ and $N$ tend to $\infty$.. This does not use the saddlepoint method,
but relies on the availability of a generating function formula for Laguerre polynomials. The proof is done in \textbf{A7}.\\
A simple modification to \citet{imj} would give the same result: in the display preceding (6.13) there, we could
write
$$
h(t)=\sum_{k=0}^{\infty}c_k t^k = 2^{\alpha/2}\Gamma(\alpha/2)(1+t)(1-t^2)^{-(\alpha/2+1)}
$$
and expand $(1-t^2)^{-(\alpha/2+1)}$. Multiplying by $1+t$ has a very simple effect on the series, and so $c_k$ is
known explicitly.

In \textbf{A8}, we show how to check that the conditions required for Soshnikov's results to hold are indeed met.
They are straightforward consequences of the analysis we will carry below.

Since the real case is derived from the complex one after analyzing a few technical points, we verify these in the
appendices and present here the study of the complex case. We now turn to the main problem we solve in this note:
showing (\ref{whattoshow1}) and (\ref{whattoshow2}) under our set of assumptions.

\section{Complex case: study of asymptotics} In this section, we work on the problem of showing
pointwise convergence and uniform boundedness, setting the problem in a way similar to section 5 of \citet{imj}.
We recall his notations, slightly modified to avoid confusions: $N_+=N+1/2$, $n_+=n+1/2$,
$z=\mu_N+\sigma_N s$, with
$\mu_N=(\sqrt{(N+\alpha)_+}+\sqrt{N_+})^{1/2}$ and $\sigma_N=(\sqrt{(N+\alpha)_+}+\sqrt{N_+})(1/\sqrt{N_+}+1/\sqrt{(N+\alpha)_+})^{1/3}$. For
reasons that will be transparent later on, our aim is to show that
\begin{equation}\label{pointwisecv}
F_N(z)=(-1)^N \sigma_N^{-1/2}\sqrt{N!/n!}\,z^{(\alpha+1)/2}\me^{-z/2}L_N^{\alpha_N}(z) \tendsto \Ai (s), \;
\forall s \in \mathbb{R} \; ,
\end{equation}
and
\begin{equation}\label{uniformcv}
F_N(z)=\gO(\me^{-s})\; \text{uniformly in } \; [s_0,\infty), s_0 \in \mathbb{R} \; .
\end{equation}

The scaling is slightly different from the original proof: $N^{-1/6}$ has been replaced by $\sN^{-1/2}$. As in
\citet{imj}, we focus on $w_N(z)=z^{(\alpha+1)/2}\me^{-z/2}L_N^{\alpha}(z)$, which satisfies
$$
\frac{d^2w}{dz^2}=\left(\frac{1}{4}-\frac{\kappa}{z}+\frac{\lambda^2-1/4}{z^2}\right)w \; ,
$$
where $\kappa=N+(\alpha+1)/2$ and $\lambda=\alpha/2$. Remark that under $\AB \definition n,N,n/N \tendsto \infty$,
$\kappa \sim \lambda$. Our strategy is to reformulate the problem in terms of so-called Whittaker functions,
denoted $W_{k,m}$,  and to use the extensive available studies of these functions to show (\ref{whattoshow1}) and
(\ref{whattoshow2}). \citet{temme}, formula (3.1) p.117 shows that
$$
w_N(z)=\frac{(-1)^N}{N!}W_{\kappa,\lambda}(z) \; .
$$

From now on, we will closely follow \citet{olver80}. Let us remark that
$$
F_N(z)=\sN^{-1/2}\frac{1}{\sqrt{n!N!}}\wkl (z) \;.
$$
We fix $s_0 \in \mathbb{R}$, and we work only with $z=\mu_N+\sN s$, where $s\geq s_0$.
\paragraph{Preliminaries} Following \citet{olver80}, we introduce $l=\kappa/\lambda$, $\beta=\sqrt{2(l-1)}$, and
the turning points $x_1=2l-2\sqrt{l^2-1}$, $x_2=2l+2\sqrt{l^2-1}$, after the rescaling $x=z/\lambda$. We remark
that the two turning points coalesce at 2 under the hypothesis $\AB$. In the new variable $x$, we have
$$
\frac{d^2W}{dx^2}=\left(\lambda^2 g(x) - \frac{1}{4x^2}\right)W \;,
$$
where $g(x)=\frac{(x-x_1)(x-x_2)}{4x^2}$. Using the ideas explained in \citet{imj}, we shall be - eventually -
interested in the asymptotics for $z=\mu_N+\sigma_N s$, or $x=z/\lambda=x_2+\sigma_N s/\lambda$ of $F_N(z)$. Let
us now define an auxiliary variable $\upsilon$ by
\begin{align*}
\int_{\beta}^{\upsilon}(\tau^2-\beta^2)^{1/2}d\tau &= \int_{x_2}^x g^{1/2}(t)dt \hspace{1cm} \text{if } x_2 \leq x
< \infty \;,
\\
\int_{-\beta}^{\upsilon}(\beta^2-\tau^2)^{1/2}d\tau &= \int_{x_1}^x (-g)^{1/2}(t)dt \hspace{1cm} \text{if } x_1
\leq x \leq x_2 \;.
\end{align*}

We limit $x$ to this range because of the technically important following point: $\sN/\lambda$ tends to zero
faster than $x_2-x_1$ does, and so, when $s$ is bounded below, $x$ will stay in the range $(x_1,\infty)$ for all
$N$ greater than a certain $N_0$. This is shown in \textbf{A2}, along with the closely related fact that we can
focus on
$\upsilon\geq 0$. Our analysis is based on section 3 of \citet{olver80}, where he builds on \citet{olver75}, in which he
expands Whittaker functions in terms of parabolic cylinder functions. The condition $\ups\geq 0$ is critical,
since Olver's expansions depend on the sign of $\upsilon$. Therefore, \textbf{A2} entitles us to focus on only one
specific form of these. From (3.10) p.219 in \citet{olver80}, one has
\begin{align*}
\wkl (\lambda x)=
(2\lambda)^{1/4}&\{\lambda(2+\beta^2/2)/e\}^{\lambda(1+\beta^2/4)}
%\\
&\times\left(\frac{\ups^2-\beta^2}{x^2-4lx+4}\right)^{1/4}x^{1/2}
\{U(-\frac{1}{2}\lambda\beta^2,\ups\sqrt{2\lambda})+\eps_1(\lambda^2,\beta^2,\ups)\} \; ,
\end{align*}
where, if $\mE$ and $\mM$ are the weight and modulus functions associated with $U$ in \citet{olver75} (p.156), we
have, according to \citet{olver80} (3.11) p.219,
\begin{equation}\label{error}
\eps_1(\lambda^2,\beta^2,\ups)=\mE^{-1}(-\frac{1}{2}\lambda\beta^2,\ups\sqrt{2\lambda})
\mM(-\frac{1}{2}\lambda\beta^2,\ups\sqrt{2\lambda})\;\mathrm{O}(\lambda^{-2/3})
\end{equation}
\textbf{uniformly} with respect to $\beta\in [0,B]$ and $\ups \in [0,\infty)$, $B$ being an arbitrary positive
constant. We recall that the main relationship between $U$, $\mE$ and $\mM$ : for $b \leq 0$ and $x\geq 0$,
$|U(b,x)|\leq \mE^{-1}(b,x)\mM(b,x)$.

 We now show that we have uniform boundedness on $[s_0,\infty)$. The pointwise convergence result will
be a straightforward consequence of the arguments we need to develop to solve this first problem.
\subsection{Uniform Boundedness}
Following up on the previous displays, if $n,N$ are large enough so that $\ups\geq 0$, we have
\begin{multline}\label{Wbound}
\left|\wkl (\lambda x)\right|\leq (2\lambda)^{1/4}\{\lambda(2+\beta^2/2)/e\}^{\lambda(1+\beta^2/4)}
\times\left(\frac{\ups^2-\beta^2}{x^2-4lx+4}\right)^{1/4}x^{1/2}\mM\mEi (1+\mathrm{O}(\lambda^{-2/3})) \;,
\end{multline}
where we omitted the argument $(-\frac{1}{2}\lambda\beta^2,\ups\sqrt{2\lambda})$ for readability purposes. Our
plan is now to transform this upper bound into a somewhat similar one, involving the modulus and weight function
associated with the Airy function, which have the advantage of having only one parameter and known asymptotics.

To carry out this program, we need to split the investigation into two parts: first $s\geq 0$ or $\ups\geq \beta$.
This will allow us to find an $s_1 \geq 0$ such that $F_N(z)=\gO(\me^{-s})$ on $[2s_1,\infty)$. In the second
part, we will just have to consider the case $s\in [s_0,2s_1]$, and show that $F_N$ is merely uniformly bounded on
this interval.
\subsubsection{Case $\mathbf{ s\geq 0} $ }
In order to use the results linking parabolic cylinder functions and the Airy function  (proved in \citet{olver59}
and cited in \citet{olver75}), let us define yet another auxiliary variable, $\eta$, by
$$
\frac{2}{3}\eta^{3/2}\beta^2=\int_{x_2}^xg^{1/2}(t)dt \;.
$$

Then, if we call $\cE$ and $\cM$ the weight and modulus functions associated with the Airy function, we have, as shown
in \textbf{A3}:
\begin{align*}
\mEi(-\frac{1}{2}\lambda\beta^2,\ups\sqrt{2\lambda}) &\leq
\cEi(\lambda^{2/3}\beta^{4/3}\eta)(1+\gO((\lambda\beta^2)^{-1})) \; ,\\
\mM(-\frac{1}{2}\lambda\beta^2,\ups\sqrt{2\lambda})&\leq
\frac{\sqrt{2}\pi^{1/4}\left(\Gamma((1+\lambda\beta^2)/2)\right)^{1/2}
\beta^{1/2}}{(\lambda\beta^2)^{1/12}}\\
&\times\left(\frac{\eta}{\ups^2-\beta^2}\right)^{1/4}\cM(\lambda^{2/3}\beta^{4/3}\eta)
\left(1+\gO((\lambda\beta^2)^{-1})\right) \; .
\end{align*}
Whence, if we call $\theta\triangleq \lambda^{2/3}\beta^{4/3}\eta$,
$$
|F_N(\lambda x)|\leq K_{n,N}x^{1/2}(\eta/(x^2-4lx+4))^{1/4}\cEi(\theta)\cM(\theta)
\left(1+\gO((\lambda\beta^2)^{-1}\vee \lambda^{-2/3})\right) \; .
$$
In \textbf{A4}, we show that $K_{n,N}\sim 2^{2/3}(N/n)^{1/4}$ under $\AB$. From now on, $\Delta$ will denote a
generic constant; its value may change from display to display. As long as $x\geq x_2$, or
$s\geq 0$, we have
$$
|F_N(\lambda x)|\leq \Delta (N/n)^{1/4} x^{1/2}(\eta/(x^2-4lx+4))^{1/4}\cEi(\theta)\cM(\theta)) \; .
$$
Now using the fact that (see \citet{olver}, chap. 11) $x^{1/4}\cM(x)\leq \Delta$, $\cEi(x)\leq \Delta
\exp(-2x^{3/2}/3)$ for $x\geq 0$ and $\lambda\beta^2=2N+1$, we get the new inequality
$$
|F_N(\lambda x)|\leq \Delta
\left(\frac{N}{n}\right)^{1/4}N^{-1/6}\left(\frac{x^2}{x^2-4lx+4}\right)^{1/4}\exp(-(2\theta^{3/2})/3) \; .
$$
In \textbf{A5.1}, we show that there exists $s_1$ such that if $s\geq 2 s_1$, $(2\theta^{3/2})/3 \geq s$. Also, as
shown in \textbf{A6.1}, if $s\geq 0$, $g$ is positive and increasing in $x$ (or, equivalently, in $s$). Since the
rational function of $x$ appearing in the previous display is just $(4g(x))^{-1/4}$, we can bound it by its
value at $x(2s_1)$ on $[2s_1,\infty)$. In \textbf{A6.2}, we show that, at $s$ fixed, under $\AB$, we have
$4g(x)\sim \beta\sN s/\lambda$, and using the equivalents mentioned in \textbf{A1}, we have $\sN\beta/\lambda\sim
4 N^{1/3}/n$, from which we conclude that
$$
\left(\frac{N}{n}\right)^{1/4}N^{-1/6}(4g(2s_1))^{-1/4}\sim N^{1/12}n^{-1/4}(8s_1N^{1/3}/n)^{-1/4}\sim
(8s_1)^{-1/4} \; .
$$
Therefore, if $N$ is large enough,
$$
\forall s\in[2s_1,+\infty)\;\;|F_N(\lambda x)| \leq \Delta \exp(-s)
$$
\subsubsection{Case $s\in[s_0,2s_1]$}\label{casesinterval}
Our aim now is just to show that $F_N$ as a function of $s$ is bounded on this interval; from this we shall
immediately have that $F_N=\gO(\exp(-s))$ on this interval, and we will have a proof of (\ref{uniformcv}).\\
This part is comparatively simpler: we use equation (\ref{Wbound}), in which we have $\mEi\leq 1$, by definition
(\citet{olver75}, p.156, (5.22)). Now using the display between (6.12) and (6.13) p.159 of the same article, we
have for $\lb\geq 1$ and $\ups\geq 0$,
$$
\frac{\mM(-\lambda\beta^2/2,\ups\sqrt{2\lambda})}{\left(\Gamma((1+\lambda\beta^2)/2)\right)^{1/2}}\leq
\frac{\Delta\beta^{1/2}}{(\lambda\beta^2)^{1/12}}\left(\frac{\eta}{\ups^2-\beta^2}\right)^{1/4} \; .
$$
Hence,
$$
|F_N(\lambda x)| \leq K_{n,N}\Delta\left(\frac{\eta}{x^2-4lx+4}\right)^{1/4}x^{1/2} \; .
$$
However on this interval, $x\tendsto 2$, by \textbf{A5.2}
$\eta=(\lambda\beta^2)^{-2/3}s+\lo((\lambda\beta^2)^{-2/3})$, and by \textbf{A6.2}
$(x^2-4lx+4)=4s\sN\beta(1+\lo(1))/\lambda$. Therefore,
$$
\frac{\eta}{x^2-4lx+4}\sim 2^{-2/3}2^{-4} n/N
$$
on the whole interval, and, because of the asymptotic estimate of $K_{n,N}$ given in \textbf{A4}, $F_N$ is bounded
uniformly in $N$ on the interval $[s_0,2s_1]$.

 We can thus conclude that
$$
\forall s_0, \; \; \exists N_0(s_0) \; N>N_0(s_0) ,\; \; \; F_N(s)=\gO_{s_0}(\me^{-s}) \; \; \text{ on }
[s_0,\infty) \; .
%, \; \; \text{\textbf{uniformly} in } N
$$
\subsection{Pointwise convergence}
Having studied in detail the uniform boundedness of $F_N$ makes the pointwise convergence problem easier. First,
since we bounded above $F_N$ in terms of $\mM$ and $\mEi$, equation (\ref{error}) shows that
$\eps_1=\gO(\lambda^{-2/3}\me^{-s})$ on $[s_0,\infty)$. So for fixed $s$, it tends to zero as $N$ gets large. The pointwise
limit of $F_N$ will be the pointwise limit of the parabolic cylinder function part of the expansion. We call this
part $\wp F_N$, for ``principal part".

Using the relationship between $U$ and $\Ai$ that we mention in \textbf{A3}, we have, with
$\theta=(\lambda\beta^2)^{2/3}\eta$,
$$
\wp F_N(\lambda x)=K_{n,N}x^{1/2}\left(\frac{\eta}{x^2-4lx+4}\right)^{1/4}(\Ai(\theta)+\cEi(\theta)\cM(\theta) \,
\gO((\lambda\beta^2)^{-1})) \; .
$$
Since $x\tendsto 2$, $K_{n,N}\sim 2^{2/3}(N/n)^{1/4}$ and given the estimate we just mentioned for the ratio
$\eta/(x^2-4lx+4)$, we have
$$
K_{n,N}\, x^{1/2}\left(\frac{\eta}{x^2-4lx+4}\right)^{1/4} \sim 1 \; .
$$
In other respects, we show in \textbf{A5.2} that $\theta\tendsto s$ under $\AB$. Finally, $\cEi$ and $\cM$ are
bounded on $\mathbb{R}$, as shown in 11.2 (pp.394-397) of \citet{olver}. Hence $\cEi(\theta)\cM(\theta) \,
(\lambda\beta^2)^{-1} \tendsto 0 $ under $\AB$, and we can conclude that $\wp F_N(\lambda x)\tendsto \Ai(s)$;
combining all the elements gives
$$
\forall s \in \mathbb{R}, \; \; \; \;F_N(\lambda x)\tendsto \Ai(s) \; \text{ under }\AB \; .
$$

\subsection{Asymptotics for $\mathbf{\phi_{\tau}} \text{ and } \mathbf{\psi_{\tau}}$}

So far we have shown that $F_N(z)=(-1)^N
\sN^{-1/2}\sqrt{z}\phi_N(z)\tendsto \Ai(s)$, and that $\me^s F_N$ was bounded when $N>N_0$ and $s \geq s_0$.\\
Our aim is to show (\ref{whattoshow1}) and (\ref{whattoshow2}). Let us write, as in \citet{imj},
$$
\phi_{\tau}=\phi_{I,N}+\phi_{I\!I,N} \; ,
$$
where
$$
\phi_{I,N}(z)=(-1)^N\sN\sqrt{a_N n}\phi_N(z)/(\sqrt{2} z)= F_N(z)d_N (z/\mu_N)^{-3/2}\;.
$$
\paragraph{Study of $\phi_{I,N}$}
In the previous display, we have $d_N=(\sN/\mu_N)^{3/2}\sqrt{a_N n/2}$. As $\sN\sim n^{1/2}N^{-1/6}$ and
$\mu_N\sim n$, $(\sN/\mu_N)^{3/2}\sim n^{-3/4}N^{-1/4}$. Since $a_N=\sqrt{Nn}$, $a_Nn=n^{3/2}N^{1/2}$, and
therefore $d_N\tendsto 1/\sqrt{2}$. But when $s$ is fixed,  $z/\mu_N \tendsto 1$, so it follows that
$$
\text{Under }\AB,\hspace{.5cm} \phi_{I,N}(\mu_N+\sN s)\tendsto \frac{\Ai(s)}{\sqrt{2}} \;.
$$
To bound $\phi_{I,N}$ for $N>N_0$ and $s\geq s_0$, we use, as in \citet{imj}, the uniform bound for $F_N$ and
$(z/\mu_N)^{-3/2}\leq \exp(-3\sN s/(2\mu_N))$, if $s\geq 0$. If $s\leq 0$,  we have $(z/\mu_N)^{-3/2}\leq
(1+s_0\sN/\mu_N)^{-3/2}$, and since this converges to 1 under $\AB$, it is bounded if $N$ is large enough. So we
have shown that,
$$
\phi_{I,N}(\mu_N+s\sN)
\left\{
\begin{array}{clc} \tendsto& 2^{-1/2}\Ai(s)\;,  \; \; N\tendsto \infty \;, \\
\leq & M \me^{-s} \; \; \text{ on } [s_0,\infty) \text{ if } N>N_0(s_0) \; .
\end{array}
\right.
$$
\paragraph{Study of $\phi_{I\!I,N}$}
We use once again the same approach as in \citet{imj}. We have
$$
\phi_{I\!I,N}=u_N v_{N-1} \phi_{I,N-1} \; ,
$$
where $u_N=(\sN/\sigma_{N-1})\sqrt{a_N/a_{N-1}}$ and $v_N=(N/n)^{1/2}$, and $n_{N-1}$ appearing in $\sigma_{N-1}$
is $n_N-1$ (for $\phi_{N-1}$ is defined in terms of  $L_{N-1}^{\alpha_N}$  and we should therefore have the same
$\alpha=n-N=(n-1)-(N-1)$). Remark that under $\AB$, $v_N \tendsto 0 $ and $u_N \tendsto 1$. \\
Define $s'$ by $\mu_N+\sN s=\mu_{N-1}+\sigma_{N-1}s'$. From
$$
s'=\frac{\mu_N-\mu_{N-1}}{\sigma_{N-1}}+\frac{\sN}{\sigma_{N-1}}s \; ,
$$
we deduce that  $s'\geq s/2 $ on $[0,\infty)$, if $N$ is large enough: as a matter of fact, under $\AB$, $\mu_N
-\mu_{N-1 }=\gO(\sqrt{n/N})$, $\sigma_{N} \sim n^{1/2}N^{-1/6}$, and $\sN/\sigma_{N-1}\tendsto 1$, so it is larger
than
$1/2$ when $N$ is large enough. To summarize, we just showed that
$$
\phi_{I\!I,N}(\mu_N+s\sN)\leq M v_N \me^{-s/2} \; \; \text{for } s\in [0,\infty) \; ,
$$
by applying the bound we got for $\phi_{I,N}$ to $\phi_{I,N-1}$ and $s'$ as the dummy variable. Here, we are
implicitly using the fact that since $n/N\tendsto \infty$, $(n-1)/(N-1)$ does too, and we can apply all the
results we derived before. On the other hand, when $s \in[s_0,0]$, we can use the fact that
$(\mu_N-\mu_{N-1})\geq 0$ and
$\sN/\sigma_{N-1} \leq 2$ to show that $s'\geq 2 s $ and hence
$$
\phi_{I\!I,N}(\mu_N+s\sN)\leq M v_N \me^{-2s}\leq M' v_N \me^{-s/2} \; \; \text{for } s\in [s_0,0] \; .
$$
The conclusion is therefore that
$$
\phi_{I\!I,N}(\mu_N+s\sN) \left\{\begin{array}{clc} \tendsto& 0 \;, \; \; N\tendsto \infty ,\\
\leq & \Delta \me^{-s/2} \; \; \text{ on } [s_0,\infty)\;, \text{ if } N>N_0(s_0) \;.
\end{array}
\right.
$$
Hence we have shown that (\ref{whattoshow1}) and (\ref{whattoshow2}) held for $\phi_{\tau}$. The analysis for
$\psi_{\tau}$ is similar.

\section{Appendices}
This section is devoted to giving background information needed to understand the problem and make the paper
relatively self-contained. We also establish many of the properties needed in the course of the proofs of
equations
(\ref{whattoshow1}) and (\ref{whattoshow2}) here.\\
Before we start, let us mention a notation issue: $\alpha$ changes value depending on whether we treat the complex
case or the real one. For the complex case $\alpha+N=n$, whereas for the real one $\alpha+N=n-1$. We frequently
replace $\alpha+N$ by $n$ in what follows; this is because the proof of equations (\ref{whattoshow1}) and
(\ref{whattoshow2}) is done in the complex case and applies to the real one by just changing $n$ into $n-1$
everywhere. When dealing with problems which are real case specific, we keep the notation
$N+\alpha$. The definition of $\mN$ and $\sN$ are also given in terms of $N+\alpha$ to highlight the adjustments
needed when dealing with the real or the complex case.
\subsection*{A0: Tracy-Widom distributions}
We recall here the definition of the Tracy-Widom distributions. We split the description according to whether
the entries of the matrix we are considering are real or complex.\\
We first need to introduce the function $q$, defined as
$$
\left\{
\begin{array}{l}
q''(x)=xq(x)+2q^3(x) \;,\\
q(x)\sim \Ai(x) \; \; \; \text{ as } x\tendsto \infty \;.
\end{array}
\right.
$$

\textbf{$\bullet$ Complex Case} The Tracy-Widom distribution appearing in the complex case, $W_2$, has cumulative
distribution function $F_2$ given by
$$
F_2(s)=\exp\left(-\int_s^{\infty}(x-s)q^2(x) dx\right)\;.
$$
The joint distribution is slightly more involved to define. Following \citet{sosh}, we do it through its
$k$-point correlation functions, using its determinantal point process character (see e.g \citet{soshDet}).\\
Let us first call $\bar{S}$ be the Airy operator. Its kernel is
$$
\bar{S}(x,y)=\frac{\Ai(x)\Ai'(y)-\Ai(y)\Ai'(x)}{x-y}=\int_0^{\infty}\Ai(x+u)\Ai(y+u) du \;.
$$
In the complex case, the $k$-point correlation functions have the property that
$$
\rho_k(x_1,\ldots,x_k)=\det_{1\leq i,j\leq k} \bar{S}(x_i,x_j) \;.
$$

\textbf{$\bullet$ Real Case} The real counterpart of $W_2$, which is called  $W_1$, has cdf $F_1$ with
$$
F_1(s)=\exp\left(-\frac{1}{2}\int_s^{\infty}q(x)+(x-s)q^2(x) dx\right) \;.
$$
The $k$-point correlation functions satisfy
$$
\rho_k(x_1,\ldots,x_k)=\left(\det_{1\leq i,j\leq k} K(x_i,x_j)\right)^{1/2} \;,
$$
where the $2\times 2$ matrix kernel of $K$ has entries (see \citet{sosh}, eq (2.18) to (2.21))
\begin{align*}
K_{1,1}(x,y)&=\bar{S}(x,y)+\frac{1}{2}\Ai(x)\int_{-\infty}^y \Ai(u) du  \; ,\\
K_{2,2}(x,y)&=K_{1,1}(y,x)\; ,\\
%\bar{S}(y,x)+\frac{1}{2}\Ai(y)\int_{-\infty}^x \Ai(u) du \\
K_{1,2}(x,y)&=-\frac{1}{2}\Ai(x)\Ai(y)-\frac{\partial}{\partial y}\bar{S}(x,y)\;,\\
K_{2,1}(x,y)&=-\int_0^{\infty}dt\left(\int_{x+t}^{\infty}\Ai(v) dv \right)\Ai(y+t)
-\eps(x-y)+\frac{1}{2}\int_y^x\Ai(u)du+\frac{1}{2}\int_x^{\infty}\Ai(u)du\int_{-\infty}^y\Ai(v)dv \; .
\end{align*}

\subsection*{A1: Asymptotic behavior of some simple functions}
In this appendix, we present some basic facts and identities that we used throughout the proof. \\
We will make repeated use of the following observations: since
$\sN=(\sqrt{(N+\alpha)_+}+\sqrt{N_+})(1/\sqrt{N_+}+1/\sqrt{(N+\alpha)_+})^{1/3} $ and $
\lambda =\alpha/2$, under \textbf{AB} we have
\begin{align*}
\sN&\sim n^{1/2}N^{-1/6} \;,\\
\lambda&\sim n/2 \; .
\end{align*}

We also use several times the following identities:
\begin{fact}\label{factequal}
With $\lambda=\alpha/2$, $\kappa=N+(\alpha+1)/2$, and $l=\kappa/\lambda$, $\beta=\sqrt{2(l-1)}$, we have
\begin{align*}
\lambda\beta^2&=(2N+1) \;,\\
\beta &\sim 2\sqrt{N/n} \;.
\end{align*}
\end{fact}
The first remark is simple algebra, and the second one comes from $\beta^2=2(l-1)=2(2N+1)/\alpha\sim 4N/n$ under $\AB$.
We have  the estimates:
\begin{fact}\label{factequiv} $x_2-x_1 \sim
8 \sqrt{N/n}$ and $\sigma_N/\lambda \sim 2 n^{-1/2}N^{-1/6}$ .
\end{fact}
The second one is obvious; the first one comes from the fact that $x_2-x_1 = 2\sqrt{2}\beta (l+1)^{1/2}$ as
$x_{2,1}=2l \pm 2\sqrt{l^2-1}$. Using Fact~\ref{factequal} immediately gives the claimed result. Finally, we have the
following estimates
\begin{fact}\label{factfunctions}
$\beta\sN/\lambda \sim 4 N^{1/3}/n$ and $\sN^3/(\lambda\beta^2)\sim(n/N)^{3/2}/2$ .
\end{fact}
The result directly follows from the aforementioned estimates.
\subsection*{A2: Working with $\ups\geq 0$}
Here we assume that $s\in[s_0,\infty)$. We also assume that $s<0$, for otherwise we can work with $\ups\geq
\beta>0$. From \textbf{A1}, we have $|x-x_2|=|s|\sN/\lambda \leq |s_0|\sN/\lambda \ll x_2-x_1 $ by
Fact~\ref{factequiv}. Now $\ups=0$ corresponds to $x_0\leq \bar{x}=(x_1+x_2)/2$: as a matter of fact, since
$(x_2-x)(x_1-x)$ is symmetric around $\bar{x}$ and $1/x$ is obviously larger on $[x_1,\bar{x}]$ than it is on
$[\bar{x},x_2]$, we have
$$
\int_{-\beta}^{\ups_{\bar{x}}} (\beta^2 - \tau^2) d\tau = \int_{x_1}^{\bar{x}}(-g(t))^{1/2}dt \geq
\int_{\bar{x}}^{x_2}(-g(t))^{1/2}dt \; .
$$
By symmetry, we also get
\begin{align*}
\int_{-\beta}^{0} (\beta^2 - \tau^2) d\tau = \int_{0}^{\beta} (\beta^2 - \tau^2) d\tau & = \frac{1}{2}
\int_{x_1}^{x_2}(-g(t))^{1/2}dt \\
&\leq \int_{x_1}^{\bar{x}}(-g(t))^{1/2}dt = \int_{-\beta}^{\ups_{\bar{x}}} (\beta^2 - \tau^2) d\tau \; ,
\end{align*}
and therefore, $\ups_{\bar{x}}>0$.\\
However $\bar{x}$ is always smaller than $x(s_0)$ if $N$ is large enough. So we can limit our investigations to
the case
$\ups\geq 0$.

\subsection*{A3: Relationship between $\mEi$, $\cEi$, $\mM$ and $\cM$}
We claim that if $s\geq 0$, and we define $\theta=(\lambda\beta^2)^{2/3}\eta$, the following inequalities hold
true:
\begin{align*}
\mEi(-\lb/2,\ups\sqrt{2\lambda})&\leq \cEi(\theta)(1+\gO((\lb)^{-1}))  \; ,\\
\mM(-\lb/2,\ups\sqrt{2\lambda})&\leq
\frac{(4\pi)^{1/4}}{(\lb)^{1/12}}[\Gamma((1+\lambda\beta^2)/2)]^{1/2}\beta^{1/2}\left(\frac{\eta}{\ups^2-\beta^2}\right)^{1/4}\\
& \times\cM(\theta)(1+\gO((\lb)^{-1})) \;.
\end{align*}

For the sake of simplicity we call $\Xi$ the part that precedes the sign ``$\times$" in the last inequality.\\
According to \citet{olver75}, equations (5.12) and (5.13), we have
\begin{align*}
U(-\lb/2,\ups\sqrt{2\lambda})&=\Xi\left\{\Ai(\theta)+\cM(\theta)\cEi(\theta)\,\gO((\lb)^{-1})\right\} \; ,\\
\bar{U}(-\lb/2,\ups\sqrt{2\lambda})&=\Xi\left\{\Bi(\theta)+\cM(\theta)\cEi(\theta)\,\gO((\lb)^{-1})\right\} \;.
\end{align*}
We have, if $s\geq 0$, $x\geq x_2$, so
\begin{equation} \label{vardef}
2/3\beta^2\eta^{3/2}=\int_{x_2}^x g^{1/2}(t)dt=\int_{\beta}^{\ups}(\tau^2-\beta^2)^{1/2}d\tau \;.
\end{equation}

For the Airy function, the weight and modulus functions had different definition depending on whether the argument was
bigger than the largest root, $c$, of $\Ai(z)=\Bi(z)$ or not. Likewise, the definition of $\mEi$ and $\mM$ depends on
the position of the argument with respect to the largest root of the equation $\bar{U}(b,x)=U(b,x)$, which is called
$\rho(b)$ in \citet{olver75}.
\paragraph{Where do the auxiliary variables lie when $\mathbf{s\geq 0}$?} We claim that the answer is that $\theta \geq 0 >
c$, and $\ups\sqrt{2\lambda}\geq \rho(-\lb/2)$.\\
The first part of equation (\ref{vardef}) implies that $\eta\geq 0$, so $\theta \geq c$, as $c<0$. This means that
we can use the definition $\cM^2=2\Ai\Bi$ and $\cEi\cM=2^{1/2}\Ai$. The second part implies that $\ups \geq
\beta$; therefore, $2\lambda\ups^2\geq 2\lb \geq \rho(-\lb/2)^2$, since by \citet{olver75}, equation (5.21),
$\rho(b)\leq 2(-b)^{1/2}$ when $b\tendsto -\infty$. This means that we have similar relationships between $\mEi$,
$\mM$, $U$, and
$\bar{U}$, to the one we had in the Airy case, $\bar{U}$ playing the role of $\Bi$, and $U$ playing the role of $\Ai$.
\paragraph{Consequences of their positions} The interesting consequence of the previous remarks is that we can
write, if $N$ is large enough, for all $s\geq 0$
$$
\mathbf{E}^{-2}(-\lb/2,\ups\sqrt{2\lambda})=\frac{U}{\bar{U}}=\frac{\cM(\theta)\cEi(\theta)}{\cM(\theta){\cal
E}(\theta)} \frac{2^{-1/2}+\gO((\lb)^{-1})}{2^{-1/2}+\gO((\lb)^{-1})} \;.
$$
In other words, we just proved that $\exists N_0$ such that $N>N_0$ implies, $\forall \, s\geq 0$
$$
\mEi(-\lb/2,\ups\sqrt{2\lambda}) \leq \cEi(\theta)(1+\gO((\lb)^{-1}) \;.
$$
By the same arguments, we derive that
$$
\mM(-\lb/2,\ups\sqrt{2\lambda})\leq \Xi \cM(\theta)(1+\gO((\lb)^{-1})\;.
$$
\subsection*{A4: Asymptotic behavior of $K_{n,N}$}
The aim here is to show that
$$
K_{n,N}\sim 2^{2/3} (N/n)^{1/4} \;.
$$
$K_{n,N}$ has the following expression:
$$
K_{n,N}=\frac{(2\lambda)^{1/4}\{\lambda(2+1/2\beta^2)/e\}^{\lambda(1+\beta^2/4)}\sqrt{2}\pi^{1/4}[\Gamma((1+\lambda\beta^2)/2)]^{1/2}\beta^{1/2}}{(\lambda\beta^2)^{1/12}\sqrt{n!N!\sN}}
\; .
$$
Since $\lambda\beta^2=(2N+1)$, $\Gamma((1+\lambda\beta^2)/2)=\Gamma(N+1)=N!$ . \\
In other respects, let $A_n=\{\lambda(2+1/2\beta^2)/e\}^{\lambda(1+\beta^2/4)}/\sqrt{n!}$ . Note that
$2\lambda+\lambda\beta^2/2=n-N+(2N+1)/2=n+1/2=n_+$. So $A_n=(n_+/e)^{n_+/2}/\sqrt{n!}$. Using Stirling's formula,
we get that $A_n\sim (n_+/n)^{n/2}(n_+/n)^{1/4}(2\pi \me)^{-1/4} \sim (2\pi)^{-1/4}$.\\
Now rewriting
$$
K_{n,N}=\frac{A_n(\lambda\beta^2)^{1/4}(8\pi)^{1/4}}{(\lambda\beta^2)^{1/12}\sqrt{\sN}} \;,
$$
we get that $K_{n,N}\sim 2^{2/3} (N/n)^{1/4}$, from using $A_n(8\pi)^{1/4}\sim \sqrt{2}$ and the second estimates
of Fact \ref{factfunctions} in \textbf{A1}.
\subsection*{A5: Asymptotic properties of $\mathbf{\eta}$}
This appendix is divided into two parts. We first show that there exists $s_1$ such that, if $s\geq 2s_1$,
\begin{gather}
\frac{2}{3}\lb \eta^{3/2} \geq s \;.\tag{P1}
\end{gather}
Then we shall show:
\begin{gather}
\text{uniformly in }s \in [a,b], \; \; \;(2N+1)^{2/3}\eta=s+\lo(1) \;.\tag{P2}
\end{gather}
\subsubsection*{A5.1: Proof of P1}
This is the argument that was used in \textbf{A8} of \citet{imj}. We repeat it for the sake of completeness. \\
Let us first suppose that $s$ is given. Since $g(x)=(x-x_1)(x-x_2)/(4x^2)$, we have
$$
\sN^2g(x) = s\frac{\sigma_N^3}{\lambda}\frac{(x_2-x_1)+s\sN/\lambda}{4(x_2+s\sN/\lambda)^2} \sim
s\frac{\sN^32\sqrt{2}\beta(l+1)^{1/2}}{16 \lambda}\sim s\frac{\beta\sN^3}{4\lambda}\;,
$$
the first equivalent coming from the fact that when $s$ is fixed, $x_2-x_1 \gg s\sN/\lambda$, and $x_2\tendsto 2$.
The second is just $l\tendsto 1$ under $\AB$. Now using the first point of Fact \ref{factfunctions} in
\textbf{A1}, together with $\sN^2\sim n N^{-1/3}$, we get that $(\beta\sN^3)/(4\lambda)\tendsto 1$. So at $s$
fixed,
$$
\sN^2 g(x)\tendsto s \;.
$$
Having this information let us now pick $s_1=8$. If $N$ is large enough, we have $\sN^2 g(x(s_1))\geq s_1/2=4$.
For all (fixed) $N$  $g$ is an increasing function of $s$. Therefore for the same $N$ we will have
$$
\forall \, s\geq s_1 \; \sN^2 g(x) \geq \sN^2 g(x(s_1))\geq s_1/2=4 \; ,
$$
and hence, since $s\geq s_1\geq 0$, $g$ is positive and we have $g^{1/2}(s)\geq 2/\sN$. Therefore,
$$
\frac{2}{3}\lb \eta^{2/3}=\lambda \int_{x_2}^xg^{1/2}(t)dt \geq \int_{x(s_1)}^xg^{1/2}(t)dt \geq
\frac{2\lambda}{\sN}\frac{\sN}{\lambda}(s-s_1)=2(s-s_1)\;.
$$
Consequently, if $s\geq 2s_1$, we have (P1).
\subsubsection*{A5.2: Proof of P2}
Without loss of generality, we can suppose that $a$ and $b$ have the same sign, and $a\geq 0$. (If it is not the
case, we can split $[a,b]=[a,0]\bigcup[0,b]$, apply the reasoning on each of these, and get the claimed result for
the original interval.)\\
The idea is that on $[a,b]$, we have
$$
\frac{(x-x_2)(x_2+b\sN/\lambda -x_1)}{4(x_2+a\sN/\lambda)^2} \geq g(x)\geq
\frac{(x-x_2)(x_2-x_1+a\sN/\lambda)}{4(x_2+b\sN/\lambda)^2} \; .
$$
Now on both sides, the terms which are not $(x-x_2)$ are $(x_2-x_1)(1+\lo(1))=4\beta(1+\lo(1))$, again because
$\sN/\lambda \ll \beta$. So if we integrate the square root of the previous inequality between $x_2$ and $x(s)$,
we get
\begin{align*}
2/3(s\sN/\lambda)^{3/2}2\sqrt{\beta}(1+\lo(1))/4 \geq 2/3 \eta^{3/2}\beta^2 \geq
2/3(s\sN/\lambda)^{3/2}2\sqrt{\beta}(1+\lo(1))/4 \;,
\end{align*}
 or
 $$
\frac{1}{2}s^{3/2}(\sN^3\beta/\lambda)^{1/2}(1+\lo(1)) \geq \eta^{3/2}\lambda\beta^2 \geq
\frac{1}{2}s^{3/2}(\sN^3\beta/\lambda)^{1/2}(1+\lo(1)) \;.
$$
The conclusion follows from \textbf{A1}, Fact~\ref{factfunctions}, whose first point, along with the estimate of
$\sN$ mentioned there, shows that $\sN^3\beta/\lambda\sim 4$. We note that (P2) also gives us pointwise
convergence of
$(\lb)^{2/3}\eta$ to $s$.

\subsection*{A6: Properties of $\mathbf{g}$}
We first show that $g$ is increasing - at $N$ fixed - as a function of $s$, if $s\geq 0$. Then we give an estimate of
$4x^2g(x)$ as $N\tendsto \infty$ and $s\in[a,b]$.
\subsubsection*{A6.1: $\mathbf{g}$ is increasing on $s\geq 0$}
Since $g(t)=(t-x_2)(t-x_1)/(4t^2)=(t^2-4lt+4)/(4t^2)$, we have
$$
g'(t)=\frac{l}{t^2}-\frac{2}{t^3}=\frac{lt-2}{t^3} \;.
$$
Now $lx_2=2l^2+2l\sqrt{l^2-1}\geq 2$, since $l=1+(2N+1)/\alpha\geq 1$. But $lx\geq lx_2$ when $s\geq 0$, and the
assertion is proved.
\subsubsection*{A6.2: On the asymptotic behavior of $\mathbf{4x^2g(x)}$ for $\mathbf{s\in [a,b]}$}
This estimate is motivated by the fact that in the course of the proof of the main result, we have to deal with an
expression of the form
$$
\frac{\eta}{x^2-4lx+4} \;.
$$
We already studied in detail $\eta$ as a function of $s$ and $N$. We now focus on $x^2-4lx+4$.

 Recalling that
$x^2-4lx+4=(x-x_2)(x-x_1)$ and $x=x_2+s\sN/\lambda$, we have
$$
x^2-4lx+4=s\frac{\sN}{\lambda}(x_2-x_1+s\frac{\sN}{\lambda})=s\frac{\sN}{\lambda}(x_2-x_1+\lo(\beta)) \;,
$$
because the first estimate in Fact \ref{factfunctions} shows that $\sN/\lambda = \lo(\beta)$, and since $s \in
[a,b]$, the previous statement holds true uniformly on this interval. Now $x_2-x_1\sim 4\beta$ under $\AB$, and
therefore, uniformly on $[a,b]$,
$$
x^2-4lx+4=s\frac{\sN}{\lambda}4\beta(1+\lo(1)) \; ,
$$
as was claimed in \ref{casesinterval}. Also, since $x=x_2+s\sN/\lambda$, and $x_2=2+(2l+2)^{1/2}\beta+\beta^2$,
$$
4g(x)=s\frac{\sN}{\lambda}\beta(1+\lo(1)) \;.
$$

\subsection*{A7: Limit of $c_{\phi}$}
Recall that under the notation of \citet{imj},
$$
\sqrt{2} c_{\phi}=\frac{1}{2}\sqrt{a_N}\left(\sqrt{N+\alpha}\int \xi_N
  - \sqrt{N} \int \xi_{N-1}\right)
$$
where
$\xi_k(x)=x^{\alpha/2-1}e^{(-x/2)}L_k^{(\alpha)}(x)\sqrt{\frac{k!}{(k+\alpha)!}}$.
We are interested in \\
\begin{align*}
v_{k,\alpha}&=\sqrt{k+\alpha}\int \xi_k - \sqrt{k} \int \xi_{k-1} \\
&=\sqrt{\frac{k!}{(k+\alpha-1)!}}\int_0^{\infty}x^{\alpha/2-1}e^{(-x/2)}
\left(L^{\alpha}_k(x)-L^{\alpha}_{k-1}(x)\right) dx \\
\text{(by \citet{szego} 5.1.13 p.102)} &=
  \sqrt{\frac{k!}{(k+\alpha-1)!}}\int_0^{\infty}x^{\alpha/2-1}e^{(-x/2)} L^{\alpha-1}_k(x) dx \\
&=\sqrt{\frac{k!}{(k+\alpha-1)!}} I_{k,\alpha} \;.
\end{align*}
Now using \citet{szego} 5.1.9 p.101,
$$
\sum_{k=0}^{\infty}w^k L_k^{\alpha-1}(x)=(1-w)^{-\alpha} \exp\left(-\frac{xw}{1-w}\right) \;.$$ So if
$F(\alpha)=\int_0^{\infty} \left(\sum_{k=0}^{\infty}w^k
  L_k^{\alpha-1}(x)\right)x^{\alpha/2-1}e^{(-x/2)} dx$, we have:
\begin{align*}
F(\alpha)&=(1-w)^{-\alpha}\int_0^{\infty} x^{\alpha/2-1}e^{(-x/2)}
\me^{-xw/(1-w)}dx\\
&=(1-w)^{-\alpha}\Gamma(\alpha/2)\left(\frac{2(1-w)}{1+w}\right)^{\alpha/2}\\
&=2^{\alpha/2}(1-w^2)^{-\alpha/2}\Gamma(\alpha/2) \;.
\end{align*}
Now if $x\geq 0,\hspace{.3cm} |L_n^{\alpha-1}(x)| \leq L_n^{\alpha-1}(-x)$, by 5.1.6 in \citet{szego}, and hence
$$
\left|\sum_{k=0}^{+\infty}w^k L_k^{\alpha-1}(x)\right|\leq \sum_{k=0}^{+\infty}|w|^k
L_k^{\alpha-1}(-x)=(1-|w|)^{-\alpha}\exp\left(\frac{x|w|}{1-|w|}\right) \;.
$$
Therefore, as long as $w \in(-1/3,1/3)$, we can switch orders of
summation, and get
$$
\sum_{k=0}^{\infty}w^k I_{k,\alpha}=2^{\alpha/2}(1-w^2)^{-\alpha/2}\Gamma(\alpha/2) \;.
$$
But
$(1-w)^{-\alpha/2}\Gamma(\alpha/2)=\sum_{k=0}^{+\infty}\frac{\Gamma(\alpha/2+k)}{k!}w^k$,
since the right-hand side converges without any difficulty on $(-1/3,1/3)$,
and hence
$$
\sum_{k=0}^{\infty}w^kI_{k,\alpha}=\sum_{m=0}^{\infty}\frac{2^{\alpha/2}\Gamma(\alpha/2+m)}{m!}w^{2m}\;.
$$
So we have
$$
\forall k\in 2\mathbb{N}, \hspace{.5cm} I_{k,\alpha}=\frac{2^{\alpha/2}\Gamma((\alpha+k)/2)}{(k/2)!} \;.
$$
Now
$v_{k,\alpha}=\sqrt{\frac{k!}{(k+\alpha-1)!}}I_{k,\alpha}=2^{\alpha/2}\frac{\Gamma((\alpha+k)/2)}{\sqrt{(k+\alpha-1)!}}\frac{\sqrt{k!}}{(k/2)!}$.
Since $\Gamma(z)\sim (z/e)^z \sqrt{2\pi/z}$, we have
\begin{align*}
\frac{\Gamma((\alpha+k)/2)}{\sqrt{\Gamma(k+\alpha)}}&\sim
2^{-(\alpha+k)/2}(\alpha+k)^{-1/4}(2\pi)^{1/4} \sqrt{2} \;, \\
\frac{\sqrt{k!}}{(k/2)!}&\sim 2^{k/2}(\pi k)^{-1/4}2^{1/4} \;,
\end{align*}
which in turn leads to
\begin{align*}
v_{k,\alpha}&\sim
2^{\alpha/2}(k(\alpha+k))^{-1/4}2^{-(\alpha+k)/2}2^{k/2} \sqrt{2}
\sqrt{2}\\
&\sim 2 (k(\alpha+k))^{-1/4} \\
&\sim 2/\sqrt{a_k} \;.
\end{align*}
Hence, as $N$ is even,  $\sqrt{2}c_{\phi}=v_{N,\alpha} \sqrt{a_N}/2 \tendsto 1 $.
\subsection*{A8: On \citet{sosh} Lemma 1}
For \citet{sosh} Lemma 1 to hold true in our case, we have to check two things. First that not only does $\sN
\phi(\mN+\sN s)\tendsto \Ai(s)/\sqrt{2}$, but also that this is true for the derivative:
\begin{equation}\label{eq:sosh1}
\sN^2 \phi'(\mN+\sN s) \tendsto \frac{1}{\sqrt{2}} \Ai'(s) \; .\tag{S1}
\end{equation}
We also have to verify that $\sN^2 \phi'(\mN+\sN s)$ is bounded above by $\Delta(s_0) \exp(-\Delta s)$ on
$[s_0,\infty)$, where $\Delta$ is a positive constant. We need to verify this for $\psi$ as well, but the techniques
are similar, so we will verify it only for $\phi$.

The second point that we need to check is that
\begin{equation}\label{eq:sosh2}
\int_0^{\infty}\left(\int_0^{z}\phi(u) du \; \psi(y+z)\right) dz \tendsto 0 \; \text { as } N \tendsto \infty \;
.\tag{S2}
\end{equation}
\subsubsection*{A8.1: Proof of (\ref{eq:sosh1})}
It is easy to see that all we need to work on are the properties of $g_N(s)=F_N(\mN+\sN s)$; if we can show that $\sN
F_N'(\mN+\sN s) \tendsto \Ai'(s)$, and that it is bounded by $\Delta(s_0) \me^{-\Delta s}$ on $[s_0,\infty)$,
we will be done. \\
We have very easily that
$$
-\sN F_N'(\mN+\sN s_1)=\int_{s_1}^{\infty} \sN^2 \left. \frac{d^2 F_N}{ds^2} \right|_{\mN+\sN u} du \; .
$$

So the strategy is clear: we want to show that the integrand in the right-hand side is bounded by an integrable
function and that it converges pointwise to $\Ai''(u)=u\Ai(u)$.\\
However,
$$
\left. \frac{d^2 F_N(x)}{dx^2} \right|_{\mN+\sN u} = \left[ \frac{1}{4}-\frac{\kN}{\mN+\sN
u}+\frac{\lambda^2-1/4}{(\mN+\sN u)^2}\right] F_N(\mN+\sN u) \;,
$$
and since we already know that $F_N(\mN+\sN s)\tendsto \Ai(s)$, we first need to check that, pointwise,
$$
\sN^2 \left[ \frac{1}{4}-\frac{\kN}{\mN+\sN s}+\frac{\lambda^2-1/4}{(\mN+\sN s)^2}\right] \tendsto s \;.
$$
In turn, this reduces to showing that
\begin{gather*}
\sN^2 \left[ \frac{1}{4}-\frac{\kN}{\mN}+\frac{\lambda^2-1/4}{\mN^2}\right] \tendsto 0  \;, \text{ and }\\
\frac{\sN^3}{\mN}\left[\frac{\kN}{\mN}-2 \frac{\lambda^2-1/4}{\mN^2}\right] \tendsto 1 \;.
\end{gather*}

The first result comes from the remarkable equality $\kN/\mN-\lambda^2/\mN^2=1/4$, which follows from the fact
that if we call
$x=\sqrt{N_+/(N+\alpha)_+}$, we have
$\kN/\mN=.5-x/(1+x)^2$ and
$\lambda^2/\mN^2=.25-x/(1+x)^2$. Using these estimates, we see that $\kN/\mN-2(\lambda/\mN)^2 =x/(1+x)^2 \sim
\sqrt{N_+/n_+}$, from which we conclude
that the second result holds.\\
Note that if we changed the centering and scaling (replacing $n$ by $\tilde{n}=n+\alpha$ and $N$ by
$\tilde{N}=N+\beta$), by studying the first expression in this case as a ``perturbation" of the study we just did, and
using the fact that $\mu_{\tilde{N}}-\mN=\gO(\sqrt{n/N})$, one could show that the first expression is then
$\gO(N^{-1/3})$, and so the result would hold. We also have corresponding results for the second expression. This shows
that we have some freedom in the centering and scaling we pick. It is also needed to show that
$$
\sN^2 \phi'(\mN+\sN s)\tendsto \Ai(s) \;,
$$
since in our splitting of $\phi$, the second part $\phi_{II,N}$ corresponds to parameters $(n-1,N-1)$, but is
centered and scaled using $\mN$ and $\sN$, defined with $(n,N)$.

To show that the sequence of functions we are interested in is bounded above by an integrable function, we split
$[s_0,\infty)$ into $[s_0,\sqrt{n}]$ and $[\sqrt{n},\infty)$. On the first interval, we can apply the previous results
since $\sN s/\mN$ is small compared to 1. So in particular the whole integrand will be smaller that $\Delta(s_0)
(1+|s|)^2 \exp(-s/2)$, after taking into account the properties of $F_N$. On the other hand, on $[\sqrt{n},\infty)$,
$\sN^2 \leq s^2$, and the denominators involving $s$ are bigger than $\mN$ and $\mN^2$ respectively, which gives
immediately that the integrand is less than $\Delta(s_0) s^2 \exp(-s/2)$. From this we conclude that the integrand
is less than $\Delta(s_0) \exp(-s/4)$, for instance, and that therefore the derivative we are interested in is
too.

It then follows easily that (\ref{eq:sosh1}) is true, and we also showed that the left-hand side of
(\ref{eq:sosh1}) is dominated on $[s_0,\infty)$ and for $N>N_0(s_0)$ by $\Delta(s_0)\me^{-s/4}$.
\subsubsection*{A8.2: Proof of (\ref{eq:sosh2})} The approach laid out in \citet{sosh} p.1044 works after some
modifications. We first write
$$
\int_0^{\infty}\left(\int_0^{z}\phi(u) du \; \psi(y+z)\right) dz = \int_0^{n^{5/8}}\left(\int_0^{z}\phi(u) du \;
\psi(y+z)\right) dz+ \int_{n^{5/8}}^{\infty}\left(\int_0^{z}\phi(u) du \; \psi(y+z)\right) dz \;.
$$
Then we can check, via a third order asymptotic development in $x$ of the right-hand side of equation (2.10) in
\citet{olver80}, that equation (2.18) therein is still true in our case, since, with his notations, $x_N\leq
n^{-3/8}$. Therefore, the analysis carried out after equation (3.21) of the same reference applies, and after
integration of the expansion following (3.22) adapted to our situation, we can show that
$$
\int_0^{n^{5/8}} \phi(u) du = \gO(n^{-n/16})
$$
With this estimate and this splitting of $[0,\infty)$, the rest of Soshnikov's argument holds true and therefore
(\ref{eq:sosh2}) can be verified.

\subsection*{A9: A quick look at special functions}
In this note, we mentioned three types of special functions, Airy, Whittaker, and parabolic cylinder functions. We
recall their definition in this appendix, as well as the main ideas behind some of the transformations Olver used.
To justify their introduction, let us say that they play a special role because it is possible, in the setting we
were in, to write the functions we studied as a perturbation of the differential equations these functions
satisfy.
\subsubsection*{A9.1: Airy function}
Let us consider the following second order differential equation:
%\begin{equation}\label{eq:airy}
$$
\frac{d^2w}{dx^2}=xw \;.
$$
%\end{equation}
\paragraph{General remark: Recessive solutions}
Since these functions are used to get asymptotic expansions, it makes sense to define the independent solutions with
respect to their behavior at $+\infty$. Usually, independent solutions $w_1$ and $w_2$ are sought, so that
$w_2=\lo(w_1)$ at a particular point of the (extended) real line. In our cases, it will be $\infty$. $w_2$ is called
the \emph{recessive} solution. That leaves the problem underdetermined, but with this in mind, one can then give
enough constraints so the problem is fully determined, and solve in terms of recessive and dominant solutions. For
a more precise definition of recessivity, see \citet{olver}, p.155.

In the case of the Airy function, we have for example: (from \citet{olver}, 11.1, p.392)
\begin{align*}
\Ai(x)&= \frac{1}{\pi}\int_0^{\infty}\cos(t^3/3+xt)dt \\
\Bi(x)&= \frac{1}{\pi}\int_0^{\infty}\{\exp(-t^3/3+xt)+\sin(t^3/3+xt)\}dt
\end{align*}
\subsubsection*{A9.2: Whittaker functions}
These are solution of the following differential equation
\begin{equation}\label{eq:whittaker}
\frac{d^2 W}{dx^2}=\left(\frac{1}{4}-\frac{\kappa}{z}+\frac{\lambda^2 -1/4}{z^2}\right)W \;.
\end{equation}
$\wkl$, the recessive solution at $\infty$, is obtained by requiring
$$
\wkl(x) \sim e^{-x/2}x^{\kappa} \; \text{as } x\tendsto \infty \;.
$$
The other solution is $M_{\kappa,\lambda}$, which is required to satisfy
$$
M_{\kappa,\lambda}(x) \sim x^{\lambda+1/2} \; \text{as } x\tendsto 0^+ \;.
$$
For more detail on these, see \citet{olver}, p.260, or \citet{olver80}.
\subsubsection*{A9.3: Parabolic cylinder functions}
According to \citet{olver59}, equation (2.9) p.133, \pcf\ satisfy (in the case we are interested in)
%\begin{equation}\label{eq:pcf}
$$
\frac{d^2 W}{dx^2}=\left(\frac{1}{4}x^2+a\right)W \;.
$$
%\end{equation}
$U(a,x)$ is chosen to satisfy
$$
U(a,x) \sim x^{-a-1/2}\me^{-x^2/4} \; \text{as } x\tendsto +\infty \;.
$$
On the other hand, $\bar{U}$ satisfies
$$
\bar{U}(a,x) \sim (2/\pi)^{1/2} \Gamma(1/2-a)x^{a-1/2}\me^{x^2/4} \; \text{as } x\tendsto +\infty\;.
$$
$\bar{U}$'s definition is actually fairly complicated, and can be found in \citet{olver59}, equation (2.12) or in
\citet{olver75}, section 5.1.
\subsubsection*{A9.4: On the usage of these functions}
As we mentioned earlier, these functions play a central role because it is relatively easy to transform the equations
in which we are interested into one of the three mentioned above, or a perturbation of it. Then a range of techniques
are available to study the effect of the perturbation, and one can sometimes, and obviously in the case we examine, get
asymptotic expansions in terms of the ``non-perturbed" solutions. Since these functions are quite
well known, information can be gathered about the function of original interest this way.\\
 For example, in
\citet{imj}, section 5, after the scaling $\xi=x/\kappa$, the Whittaker equation (\ref{eq:whittaker}) becomes
$$
\frac{d^2W}{d\xi^2}=\left(\kappa^2\frac{(\xi-\xi_1)(\xi-\xi_2)}{4\xi^2}-\frac{1}{4\xi^2}\right)W \;.
$$
Using the Liouville-Green transformation $\zeta(d\zeta/d\xi)^2=(\xi-\xi_1)(\xi-\xi_2)/(4\xi^2)$, with
$w=(d\zeta/d\xi)^{-1/2} W$, one has
$$
\frac{d^2w}{d\zeta^2}=\{\kappa^2 \zeta + \psi(\zeta)\} w \;.
$$
This is a perturbation of the (scaled) Airy equation, for $\Ai(\kappa^{2/3} \zeta)$ and $\Bi(\kappa^{2/3} \zeta)$ are
solutions of $d^2w/d\zeta^2=\kappa^2 \zeta w$.\\
$w$ is not $W$, but it can be related to it, and it is through this mean that Johnstone did his original analysis.
As $\psi$ is a relatively involved function of $\xi$ and $\zeta$, we do not explicit it, but just mention that the
understanding of $\psi$ is key to getting the uniform bound (\ref{whattoshow2}). For more on this, see \citet{imj}
or \citet{olver}, theorem 11.3.1 p.399.

The problem we encountered (and mentioned in \ref{errorcontrolproblem}) about the error control function is
exactly here: we could not get enough information about the behavior of $\psi$ under $\AB$, so we slightly changed
approach and turned to other studies.

In \citet{olver80}, Olver starts with equation (\ref{eq:whittaker}), where the dummy variable was $z$. Writing
$x=z/\lambda$ and $l=\kappa/\lambda$, he gets
$$
\frac{d^2 W}{dx^2}=\left(\lambda^2g(x)-\frac{1}{4x^2}\right)W \;.
$$
As he aims to expand the solution in terms of parabolic cylinder functions, he changes variables another time, by
writing
$$
W=\left(\frac{dx}{d\zeta}\right)^{1/2} w \; \; \;, \;
\left(\frac{d\zeta}{dx}\right)^2=\frac{x^2-4lx+4}{4x^2(\zeta^2-\beta^2)} \;,
$$
with $\beta=\{2(l-1)\}^{1/2}$. Hence, he gets
$$
\frac{d^2w}{d\zeta^2}=\{\kappa^2 (\zeta^2-\beta^2) + \psi(\kappa,\beta,\zeta)\} w \;,
$$
with $\psi(\kappa,\beta,\zeta)=-\dot{x}^2/(4x^2)+\dot{x}^{1/2}d^2(\dot{x}^{-1/2})/d\zeta^2$. His \citet{olver75}
is a study of this type of equations, and in particular of the control of the deviation of the solution of the
previous equation to the corresponding parabolic cylinder function. In \citet{olver80}, he studies very explicitly
the abstract estimate he gets in \citet{olver75} in the case of Whittaker functions. We use this repeatedly in our
study, as it is essential to get the crucial property (\ref{whattoshow2}).

\bibliographystyle{plainnat}
%\nocite{*}

\bibliography{research}

\end{document}